%% file: main.tex
\documentclass[reqno, a4paper, 12pt]{amsart}

\usepackage{microtype}
\usepackage[utf8]{inputenc}
\input{preamble.tex}

\newtheoremstyle{case}{}{}{}{}{}{:}{ }{}
\theoremstyle{case}
\newtheorem{case}{Case}

\begin{document}

\title[Covering
$3$-edge-coloured random graphs with monochromatic trees]{Covering
$3$-edge-coloured random graphs with monochromatic trees}

\author[Y.~Kohayakawa]{Yoshiharu Kohayakawa}
\address{Instituto de Matem\'{a}tica e Estat\'{\i}stica, Universidade de S\~{a}o
  Paulo, Brazil} 
\email{yoshi@ime.usp.br}

\author[W.~Mendon\c{c}a]{Walner Mendon\c{c}a}
\address{Instituto Nacional de Matem\'{a}tica Pura e Aplicada, Rio de Janeiro,
Brazil}
\email{walner@impa.br}

\author[G.~O.~Mota]{Guilherme Oliveira Mota}
\address{Instituto de Matem\'{a}tica e Estat\'{\i}stica, Universidade de S\~{a}o
  Paulo, Brazil} 
\email{mota@ime.usp.br}

\author[B.~Sch\"{u}lke]{Bjarne Sch\"{u}lke}
\address{Fachbereich Mathematik, Universit\"{a}t Hamburg, Hamburg,
Germany}
\email{bjarne.schuelke@uni-hamburg.de}

\thanks{\rule[-.2\baselineskip]{0pt}{\baselineskip}%
  Y.~Kohayakawa was partially supported by CNPq (311412/2018--1, 423833/2018--9) and FAPESP (2018/04876--1). 
  G.~O.~Mota was partially supported by CNPq
  (304733/2017--2, 428385/2018--4) and FAPESP (2018/04876--1). 
  W.~Mendon\c{c}a was partially supported by CAPES (88882.332408/2010--01).
  B.~Sch\"ulke was partially supported by G.I.F. Grant Agreements
  No.~I-1358--304.6/2016. 
  The cooperation of the authors was supported by a joint CAPES/DAAD
  PROBRAL project (Proj.~430/15, 57350402, 57391197).
  This study was financed in part by CAPES, Coordenação
  de Aperfeiçoamento de Pessoal de Nível Superior, Brasil,
  Finance Code 001.
  FAPESP is the S\~ao Paulo Research Foundation.  CNPq is the National
  Council for Scientific and Technological Development of Brazil.%
}

\subjclass{05C80 (primary); 05C70 (secondary)}
\keywords{Graph partitioning, monochromatic trees, random
  graphs}

\begin{abstract}
  We investigate the problem of determining how many monochromatic trees are
  necessary to cover the vertices of an edge-coloured random graph. More
  precisely, we show that for $p\gg n^{-1/6}{(\ln n)}^{1/6}$, in any
  $3$-edge-colouring of the random graph $G(n,p)$ we can find three monochromatic
  trees such that their union covers all vertices. This improves, for three
  colours, a result of Buci\'c, Kor\'andi and Sudakov.
\end{abstract}

\maketitle

\section{Introduction}

Given a graph $G$ and a positive integer $r$, let $\tc_r(G)$ denote the minimum
number $k$ such that in any $r$-edge-colouring of $G$, there are $k$
monochromatic trees $T_1,\ldots,T_k$ such that the union of their vertex sets covers $V(G)$,
i.e.,
\begin{align*}
	V(G) = V(T_1)\cup \dots \cup V(T_k).
\end{align*}
We define $\tp_r(G)$ analogously by requiring the union above to be disjoint.

It is easy to see that $\tp_2(K_n) = 1$ for all $n\geq 1$, and Erd\H{o}s,
Gy\'arf\'as~and Pyber~\cite{EGP91} proved that $\tp_3(K_n) = 2$ for all
$n\geq 1$, and conjectured that $\tp_r(K_n)=r-1$ for every $n$ and $r$. Haxell
and Kohayakawa~\cite{HaKo96} showed that $\tp_r(K_n) \leq r$ for all
sufficiently large $n \ge n_0(r)$. We remark that it is easy to see that
$\tc_r(K_n) \leq r$ (just pick any vertex $v \in V(K_n)$ and let~$T_i$, for
$i\in[r]$, be a maximal monochromatic tree of colour $i$ containing $v$), but it
is not even known whether or not $\tc_r(K_n) \leq r-1$ for every $n$ and $r$ (as
would be implied by the conjecture of Erd\H{o}s, Gy\'arf\'as~and Pyber).

Concerning general graphs instead of complete graphs, Gy\'arf\'as~\cite{Gyarfas}
noted that a well-known conjecture of Ryser on matchings and transversal sets in
hypergraphs is equivalent to the statement that for every graph $G$ and
integer~$r\geq 2$, we have $\tc_r(G)\leq (r-1)\alpha(G)$. In particular, Ryser's
conjecture, if true, would imply that $\tc_r(K_n) \leq r-1$, for every $n\geq 1$
and $r \geq 2$. Ryser's conjecture was proved in the case $r = 3$ by
Aharoni~\cite{Aharoni}, but for $r\geq 4$ very little is known. For example,
Haxell and Scott~\cite{HaSc12} proved (in the context of Ryser's original
conjecture) that there exists $\epsilon >0$ such that for $r\in \{4,5\}$, we
have $\tc_r(G)\leq (r-\epsilon)\alpha(G)$, for any graph $G$.

Bal and DeBiasio~\cite{BaDe} initiated the study of covering and partitioning
random graphs by monochromatic trees. They proved that if $p\ll
{\left(\frac{\log n}{n}\right)}^{1/r}$, then with high probability\footnote{We
  will write shortly \emph{w.h.p.} for \emph{with high probability}.} we have
$\tc_r(\gnp) \to \infty$. They conjectured that for any $r\geq 2$, this was the
correct threshold for the event $\tp_r(\gnp) \leq r$. Kohayakawa, Mota and
Schacht~\cite{KoMoSc} proved that this conjecture holds for $r=2$, while Ebsen,
Mota and Schnitzer\footnote{A description of this construction can be found
  in~\cite{KoMoSc}.} showed that it does not hold for more than two colours.

Buci\'c, Kor\'andi and Sudakov~\cite{BuKoSu} proved that if $p \ll
{\left(\frac{\log n}{n}\right)}^{\sqrt r/2^{r-2}}$, then w.h.p.\ we have
$\tc_r(\gnp) \geq r+1$, which implies that the threshold for the event
$\tc_r(G)\leq r$ is in fact significantly larger than the one conjectured by Bal
and DeBiasio when $r$ is large. Buci\'c, Kor\'andi and Sudakov also proved that
w.h.p.\ we have $\tc_r(\gnp) \leq r$ for $p \gg {\left(\frac{\log
      n}{n}\right)}^{1/2^r}$. They were also able to roughly determine the
typical behaviour of $\tc_r(G(n,p))$ in terms of the range where $p$ lies in
(see~\cite[Theorem~1.3 and Theorem~1.4]{BuKoSu}).

Considering colourings with three colours, the results from~\cite{BuKoSu} imply
that if $p \gg {\left( \frac{\log n}{n} \right)}^{1/8}$, then w.h.p.\ we have
$\tc_3(G(n,p))\leq 3$, and if ${\left( \frac{\log n}{n} \right)}^{1/6} \ll p \ll
{\left(\frac{\log n}{n} \right)}^{1/7}$, then w.h.p.\ $\tc_3(G(n,p))\leq 88$. Our
main result improves these bounds for three colours.
\begin{theorem}\label{main_res}
  If $p = p(n)$ satisfies $p \gg {\big(\frac{\log n}{n}\big)}^{1/6}$, then with high probability we
  have
  \[
    \tc_3\big( G(n,p) \big) \leq 3.
  \]
\end{theorem}

It is easy to see that if $p = 1 - \omega(n^{-1})$, then w.h.p.\ there is a
$3$-edge-colouring of~$G(n,p)$ for which three monochromatic trees are needed to
cover all vertices --- it suffices to consider three non-adjacent vertices
$x_1$, $x_2$ and $x_3$, and colour the edges incident to~$x_i$ with colour~$i$
and colour all the remaining edges with any colour. Therefore, the bound for
$tc_3(\gnp)$ in~\cref{main_res} is the best possible as long as $p$ is not too
close to $1$.

We remark that, from the example described in~\cite{KoMoSc}, we know that for $p
\ll {\left( \frac{\log{n}}{n} \right)}^{1/4}$, we have w.h.p.\ $\tc_3(G(n,p))
\geq 4$. It would be very interesting to describe the behaviour of $\tc_3(\gnp)$
when ${\big(\frac{\log n}{n}\big)}^{1/4}\ll p \ll {\big(\frac{\log
    n}{n}\big)}^{1/6}$.

This paper is organized as follows. In~\cref{sec:preliminaries} we present some
definitions and auxiliary results that we will use in the proof of
Theorem~\ref{main_res}, which is outlined in~\cref{sec:sketch}. The details of
the proof of~\cref{main_res} are given in Section~\ref{sec:proof}.

\section{Preliminaries}\label{sec:preliminaries}
Most of our notation is standard (see~\cites{Bo98,BM08,Di10}
and~\cites{Bo01,JLR00}). However, we will mention in the following few
definitions regarding hypergraphs that will play a major role in our proofs just
for completeness.

We say that a set $A$ of vertices in a hypergraph $\mathcal{H}$ is a
\emph{vertex cover} if every hyperedge of $\mathcal{H}$ contains at least one
element of $A$. The \emph{covering number} of $\mathcal{H}$, denoted by
$\tau(\mathcal{H})$, is the smallest size of a vertex cover in $\mathcal{H}$. A
\emph{matching} in $\mathcal{H}$ is a collection of disjoint hyperedges in
$\mathcal{H}$. The \emph{matching number} of $\mathcal{H}$, denoted
by~$\nu(\mathcal{H})$, is the largest size of a matching in $\mathcal{H}$. An
immediate relationship between $\tau(\mathcal{H})$ and $\nu(\mathcal{H})$ is the
inequality~$\nu(\mathcal{H}) \leq \tau(\mathcal{H}) $. If additionally
$\mathcal{H}$ is $r$-uniform, then we have~$\tau(\mathcal{H}) \leq r
\nu(\mathcal{H})$. A conjecture due to Ryser (which first appeared in the thesis
of his Ph.D. student, Henderson~\cite{Ryser}) states that for every $r$-uniform
$r$-partite hypergraph $\mathcal{H}$, we have $\tau(\mathcal{H})\leq
(r-1)\nu(\mathcal{H})$. Note that K\"{o}nig-Egerv\'{a}ry theorem corresponds to
Ryser's conjecture for $r=2$. Aharoni~\cite{Aharoni} proved that Ryser's
conjecture holds for $r=3$, but the conjecture remains open for $r\geq 4$.

Given a vertex $v$ in a 3-uniform hypergraph $\mathcal{H}$, the \emph{link
  graph} of $\mathcal{H}$ with respect to~$v$ is the graph $L_v = (V,E)$ with
vertex set $V = V(\mathcal{H})$ and edge set $E = \{xy : \{x,y,v\} \subseteq
\mathcal{H}\}$.

We will use the following theorem due to Erd\H{o}s, Gy\'arf\'as and
Pyber~\cite{EGP91} in the proof of our main result.
\begin{theorem}[Erd\H{o}s, Gy\'arf\'as and
  Pyber]\label{lemma:EGP91} For any
  $3$-edge-colouring of a complete graph $K_n$, there exists a partition of
  $V(K_n)$ into~$2$ monochromatic trees.
\end{theorem}

We will also use the following lemma, which is a simple application of
Chernoff's inequality. For a proof of the first item see~\cite[Lemma~3.8]{KMNSS}.
The second item is an immediate corollary of~\cite[Lemma~3.10]{KMNSS}.
\begin{lemma}\label{lem:gnp}
  Let~$\eps > 0$. If~$p=p(n)\gg{\left(\frac{\log{n}}{n}\right)}^{1/6}$, then
  w.h.p.\ $G\in G(n,p)$ has the following properties.
  \begin{enumerate}[label=\rmlabel]
  \item For any disjoint sets~$X,Y\subseteq V(G)$ with~$|X|,|Y|
    \gg \frac{\log{n}}{p}$, we have
    \begin{align*}
    |E_G(X,Y)| = (1\pm \eps)p|X||Y|.
    \end{align*}
  \item Every vertex $v\in V(G)$ has degree $d_G(v)=(1\pm
    \eps)pn$ and
    every set of~$i\leq 6$ vertices has $(1\pm\eps) p^i n$ common neighbours.
  \end{enumerate}
\end{lemma}

\section{A sketch of the proof}\label{sec:sketch}

In this section we will give an overview of the proof of~\cref{main_res}.
Let $G = G(n,p)$, with $p \gg {\left( \frac{\log{n}}{n} \right)}^{1/6}$, and let
$\phi: E(G) \to \{\red,~\green,~\blue\}$ be any 3-edge-colouring of~$G$. We consider an
auxiliary graph~$F$, with~$V(F)=V(G)$ and~$ij \in E(F)$ if and only if there is,
in the colouring~$\phi$, a monochromatic path in~$G$ connecting $i$ and $j$.
Then we define a 3-edge-colouring~$\phi'$ of~$F$ with~$\phi'(ij)$ being the color of any
monochromatic path in~$G$ connecting~$i$ and~$j$. Note that any covering of~$F$
with monochromatic trees with respect to the colouring~$\phi'$ corresponds to a
covering of~$G$ with monochromatic trees with respect to the colouring~$\phi$ with
the same number of trees.

Next, we consider different cases depending on the value of~$\alpha(F)$.
If~$\alpha(F)=1$, then~$F$ is a complete $3$-edge-coloured graph and by a
theorem of Erd\H{o}s, Gy\'{a}rf\'{a}s and Pyber (see~\cref{lemma:EGP91}), there
exists a partition of~$V(F)$ into~$2$ monochromatic trees. The remaining proof
now is divided into the cases $\alpha(F) \geq 3$ and $\alpha(F) = 2$.

\medskip
\noindent \textit{Case $\alpha(F) \geq 3$.} From the condition on the independence number of $G$, there
exist three vertices $r,b,g\in V(G)$ that pairwise do not have any monochromatic
path connecting them. With high probability, they have a common neighbourhood in
$G$ of size at least $np^3/2$.
 Let~$X_{rbg}$ be the largest subset of this common
neighbourhood such that for each~$i\in\{r,b,g\}$, the edges from $i$ to $X_{rbg}$ in
$G$ are all coloured with one colour. Then, since there are no monochromatic
paths between any two of $r$, $b$, $g$, we have $|X_{rbg}| \geq np^3/12$ and moreover
we may assume that all edges between~$r$ and~$X_{rbg}$ are red, all between~$b$ and~$X_{rbg}$ are
blue and those between~$g$ and~$X_{rbg}$ are green. Now we notice that all vertices that
have a neighbour in~$X_{rbg}$ are covered by the union of the spanning trees of the
red component of~$r$, the blue component of~$b$ and the green component of~$g$.

We are done in the case where every vertex has a neighbour in $X_{rbg}$, as the vertices in $X_{rbg} \cup N_G(X_{rbg})$ are covered by the red, blue and green
component containing $r$, $b$ and $g$, respectively.
Otherwise, w.h.p.\ any vertex $y\in V \setminus \left( X_{rbg} \cup N_G(X_{rbg}) \right)$ has many common neighbours with~$r$,~$g$ and~$b$ in $G$ that are
also neighbours of some vertex in $X_{rbg}$.
 An analysis of the possible colourings of the edges between $X_{rbg}$ and the common neighbourhood of the vertices $r$, $b$, $g$ and $y$ yields the following: for some $i \in \{r,g,b\}$, let us say $i=r$, every vertex $y \in X_{rbg}$
can be connected to~$r$ by a monochromatic path in colour~$\red$ or either
to~$g$ or~$b$ by a monochromatic path in the colour~$\blue$ or~$\green$,
respectively.

This already gives us that all vertices in~$G$ can be covered by~$5$
monochromatic trees, since all the vertices in~$N_G(X_{rbg})$ lie in the $\red$ component
of~$r$, or the $\green$ component of~$g$, or in the $\blue$ component of~$b$ and every
vertex in $V\setminus N_G(X_{rbg})$ lies in the $\red$ component of~$r$, in the $\blue$ component of $g$ or in the $\green$
component of $b$. By analysing the colours of edges to the common neighbourhood
of carefully chosen vertices, we are able to show that actually three of those
five trees already cover all the vertices of $G$.

\medskip
\noindent \textit{Case $\alpha(F) =2$.} Let us consider a $3$-uniform hypergraph
$\mathcal{H}$ defined as follows (this definition is inspired by a construction
of Gy\'{a}rf\'{a}s~\cite{Gyarfas}). The vertices of $\mathcal{H}$ are the
monochromatic components of~$F$ and three vertices form a hyperedge if the
corresponding three components have a vertex in common (in particular, those
three monochromatic components must be of different colours).
Hence~$\mathcal{H}$ is an~$3$-uniform $3$-partite hypergraph.

We observe that if~$A$ is a vertex cover of~$\mathcal{H}$, then the
monochromatic components associated with the vertices in~$A$ cover all the
vertices of~$G$. This implies that $\tc_{3}(G) \leq \tau(\mathcal{H})$. Also, it
is easy to see that $\nu(\mathcal{H}) \leq \alpha(F)=2$. Now, recall that
Aharoni's result~\cite{Aharoni} (which corresponds to Ryser's conjecture for
$r=3$) states that for every $3$-uniform $3$-partite hypergraph $\mathcal{H}$ we
have $\tau(\mathcal{H}) \leq 2\nu(\mathcal{H})$. Together with the previous
observation, this implies $\tc_3(G) \leq 4$. But our goal is to prove that
$\tc_3(G) \leq 3$. To this aim, we analyze the hypergraph $\mathcal{H}$ more
carefully, reducing the situation to a few possible settings of components
covering all vertices. In each of those cases, we can again analyse the possible
colouring of edges of common neighbours of specific vertices, inferring that
indeed there are~$3$ monochromatic components cover all vertices.

\section{Proof of~\cref{main_res}}\label{sec:proof}

Instead of analysing the colouring of the graph $G=G(n,p)$, it will be helpful
to analyse the following auxiliary graph.

\begin{definition}[Shortcut graph]
  Let $G$ be a graph and $\varphi$ be a $3$-edge-colouring of $G$. The
  \emph{shortcut graph} of $G$ (with respect to $\varphi$) is the graph
  $F=F(G,\phi)$ that has $V(G)$ as the vertex set and the following edge set:
  \[\{uv : u,v \in V(G) \text{ and $u$ and $v$ are connected in $G$ by a path
      monochromatic under $\phi$}\}.\]
\end{definition}

We can consider a natural edge colouring~$\phi'$ of $F(G,\phi)$ by assigning to
an edge~$uv\in E(F(G,\phi))$ the colour of any monochromatic path connecting~$u$
and~$v$ in~$G$ under the colouring~$\phi$. We will say that~$\phi'$ is an
\emph{inherited colouring} of $F(G,\phi)$. Let~$\tc(F,\phi')$ be the minimum
number of monochromatic components (under the colouring~$\phi'$) covering all
the vertices of~$F$. Note that any covering of~$F(G,\phi)$ with monochromatic
trees under~$\phi'$ corresponds to a covering of~$G$ with monochromatic trees
under the colouring~$\phi$. In particular, if we show that for every
$3$-edge-colouring $\phi$ of $G$, we have $\tc(F,\phi')\leq 3$, for every
ineherited colouring $\phi'$, then we have shown that $\tc_3(G) \leq 3$.
Therefore,~\cref{main_res} follows from the following lemma.

\begin{lemma}\label{lemma:main}
  Let $p\gg{\left(\frac{\log n}{n}\right)}^{1/6}$ and let $G =\gnp$. The
  following holds with high probability. For any $3$-edge-colouring $\phi$ of
  $G$ and any inherited colouring $\phi'$ of the shortcut graph $F = F(G,\phi)$,
  we have $\tc(F,\phi')\leq 3$.
\end{lemma}

The proof of~\cref{lemma:main} is divided into two different cases,
depending on the independence number of $F$. Subsections~\ref{sec:alpha2}
and~\ref{sec:alpha3} are devoted, respectively, to the proof of
\cref{lemma:main} when $\alpha(F)\geq 3$ and $\alpha(F)\leq 2$.

From now on, we fix $\eps>0$ and assume that~$p\gg{\left(\frac{\log
  n}{n}\right)}^{1/6}$ and $n$ is sufficiently large. Then, by~\cref{lem:gnp},
we may assume that the following holds w.h.p.:	
\begin{enumerate}
\item There is an edge between any two sets of size~$\omega\left((\log n)/p\right)$.
\item Every vertex $v\in V(G)$ has degree $d_G(v)=(1\pm \eps)pn$.
\item Every set of~$i\leq 6$ vertices has $(1\pm\eps) p^i n$ common neighbours.
\end{enumerate}

\subsection{Shortcut graphs with independence number at least three}\label{sec:alpha2}

\begin{proof}[Proof of~\cref{lemma:main} for $\alpha(F)\geq 3$]
  Since $\alpha (F) \geq 3$, there exist three vertices~$r,b,g\in V(G)$ that
  pairwise do not have any monochromatic path connecting them in~$G$. In
  particular, if $v$ is a common neighbour of $r$, $b$ and $g$ in $G$, then the
  edges $vr$, $vb$ and $vg$ have all different colours. The common neighbourhood
  of $r$, $b$ and $g$ in $G$ has size at least~$np^3/2$. Let~$X_{rbg}$ be the
  largest subset of this common neighbourhood such that for each~$i\in
  \{r,b,g\}$, the edges between $i$ and the vertices of $X_{rbg}$ are all
  coloured with the same colour in~$G$. Then~$\vert X_{rbg}\vert\geq np^3/12$.
  Without loss of generality, assume that all edges between~$r$ and the vertices
  of~$X_{rbg}$ are red, between~$b$ and the vertices of~$X_{rbg}$ are blue and
  those between~$g$ and the vertices of~$X_{rbg}$ are green.
  Let~$C_\red(r)$,~$C_\blue(b)$ and~$C_\green(g)$ be respectively the red, blue
  and green components in $G$ containing~$r$,~$g$ and~$b$.
		
  Notice that all vertices of~$F$ that have a neighbour in~$X_{rbg}$ are covered
  by~$C_\red(r)$,~$C_\blue(b)$ or~$C_\green(g)$. Therefore, the proof would be
  finished if every vertex had a neighbour in~$X_{rbg}$. If this is not the
  case, we fix an arbitrary vertex~$y\in V\setminus \left(X_{rbg} \cup
    N_G(X_{rbg}) \right)$. By our choice of~$p$, there are at least~$np^4/2$
  common neighbours of~$y$,~$r$,~$b$ and~$g$. Let~$X_{yrbg}$ be the largest
  subset of the common neighbourhood of~$y$,~$r$,~$b$ and~$g$ such that for
  each~$i\in\{r,b,g\}$, the edges between~$i$ and~$X_{yrbg}$ are all coloured
  the same. Then~$|X_{yrbg}|\geq np^4/12$. Note that since~$y\notin
  N_G(X_{rbg})$, the sets~$X_{yrbg}$ and~$X_{rbg}$ are disjoint. Furthermore,
  since~$|X_{yrbg}|,|X_{rbg}| \gg \frac{\log{n}}{p}$, we have
  \begin{align*}
    |E_G(X_{yrbg},X_{rbg})|\geq 1.
  \end{align*}

  We now analyse the colours between~$r$, $b$, $g$ and the set~$X_{yrbg}$.
  Again, since there is no monochromatic path connecting any two of~$r$,~$b$
  and~$g$, all~$i \in \{r,b,g\}$ have to connect to~$X_{yrbg}$ in different
  colours. Since~$X_{yrbg}$ is disjoint of~$X_{rbg}$, we cannot have~$r$,~$b$
  and~$g$ being simultaneously connected to~$X_{yrbg}$ by red, blue and green
  edges, respectively. Assume first that for each~$i \in \{r,b,g\}$, the edges
  between~$i$ and~$X_{yrbg}$ have different colours from the edges between~$i$
  and~$X_{rbg}$. Then let~$uv$ be an edge between~$X_{yrbg}$ and~$X_{rbg}$ and
  notice that whatever the colour of~$uv$ is, we will have a monochromatic path
  connecting two of the vertices in~$\{r,g,b\}$. Therefore, we can assume that
  for some~$i \in \{r,g,b\}$, we have that all the edges between~$i$
  and~$X_{rbg}$ and all the edges between~$i$ and~$X_{yrbg}$ coloured the same.
  Without loss of generality, we may say that such~$i$ is~$r$. In this case, the
  edges between~$b$ and~$X_{yrbg}$ are green and the edges between~$g$
  and~$X_{yrbg}$ are blue. Finally, all the edges between~$X_{yrbg}$ and~$X_{rbg}$ are
  red, otherwise we would be able to connect~$b$ and~$g$ by some monochromatic
  path. Figure~\ref{fig:rbgy} shows the colouring of the edges that we have
  analysed so far.
  
  \begin{figure}
    \centering
    \begin{tikzpicture} [scale=1, thick, auto, vertex/.style={circle, draw,
        fill=black!50, inner sep=0pt, minimum width=4pt}]

      \node [vertex, label=left:$r$](r) at (0,2) {};
      \node [vertex, label=left:$b$](b) at (0,0) {};
      \node [vertex, label=left:$g$](g) at (0,-2) {};
      \node [vertex, label=right:$y$](y) at (6,0) {};
      \node [vertex, label={[label distance=.4cm]-45:$X_{rbg}$}](j) at (4,-1) {};
      \node [vertex, label={[label distance=.2cm]+45:$X_{yrbg}$}](i) at (4,1) {};
      
      \draw (j) circle (0.6 cm);
      \draw (i) circle (0.4 cm);
      
      \draw [red!75!black] (r)--(j);
      \draw [blue!75!black] (b)--(j);
      \draw [green!75!black] (g)--(j);
      
      \draw (y)--(i);
      \draw [red!75!black](j)--(i);
      \draw [red!75!black] (r)--(i);
      \draw [green!75!black] (b)--(i);
      \draw [blue!75!black] (g)--(i);
      
    \end{tikzpicture}
    \caption{Analysis of the colouring of the edges incident on $X_{rbg}$ and on
      $X_{yrbg}$.}\label{fig:rbgy}
  \end{figure}
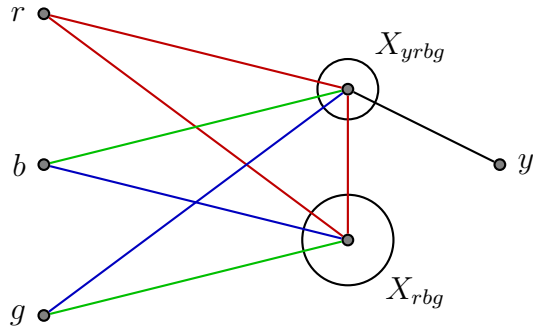
       
  Let us now consider any further vertex~$x\in V\setminus \left(X_{rbg} \cup
    N_G(X_{rbg}) \right)$ with~$x\neq y$, if such a vertex exists. We
  define~$X_{xrbg}$ analogously to~$X_{yrbg}$ and observe that the colour
  pattern from~$r$,~$b$,~$g$ to~$X_{xrbg}$ must be the same as the one
  to~$X_{yrbg}$. Indeed, if this is not the case, then a similar analysis of the
  colours of the edges between~$\{r,b,g\}$ and~$X_{xrbg}$ yields that for
  some~$i \in \{b,g\}$, we know that the edges connecting~$i$ to~$X_{xrbg}$ are
  of the same colour as the edges connecting~$i$ to~$X_{rbg}$. Without loss of
  generality, let us say that~$i$ is~$g$. Then the edges between~$b$
  and~$X_{xrbg}$ are red and the edges between~$r$ and~$X_{xrbg}$ are green,
  otherwise~$X_{xrbg}$ and~$X_{rbg}$ would not be disjoints sets.
  Figure~\ref{fig:rbgyx} shows the colouring of the edges incident to~$X_{yrbg}$
  and~$X_{xrbg}$. Since~$|X_{yrbg}|,|X_{xrbg}| \gg \frac{\log{n}}{p}$, we have
  that there is some edge~$uv$ between~$X_{yrbg}$ and~$X_{xrbg}$. But then
  however we colour~$uv$, we will get an monochromatic path connecting two
  vertices in~$\{r,b,g\}$, which is a contradiction. Thus, the colour pattern of
  edges between~$\{r,b,g\}$ and~$X_{xrbg}$ is the same as the colour pattern of
  the edges between~$\{r,b,g\}$ and~$X_{yrbg}$.

  \begin{figure}
    \centering
    \begin{tikzpicture} [scale=1, thick, auto, vertex/.style={circle, draw,
        fill=black!50, inner sep=0pt, minimum width=4pt}]

      \node [vertex, label=above:$r$](r) at (0,2) {};
      \node [vertex, label=above:$b$](b) at (0,0) {};
      \node [vertex, label=above:$g$](g) at (0,-2) {};
      \node [vertex, label=right:$y$](y) at (6,1) {};
      \node [vertex, label=right:$x$](x) at (6,-1) {};
      \node [vertex, label={[label distance=.4cm]left:$X_{rbg}$}](j) at (-2,0) {};
      \node [vertex, label={[label distance=.2cm]+45:$X_{yrbg}$}](i) at (4,1) {};
      \node [vertex, label={[label distance=.2cm]-45:$X_{xrbg}$}](k) at (4,-1) {};
      
      \draw (j) circle (0.6 cm);
      \draw (i) circle (0.4 cm);
      \draw (k) circle (0.4 cm);
      
      \draw [red!75!black] (r)--(j);
      \draw [blue!75!black] (b)--(j);
      \draw [green!75!black] (g)--(j);
      
      \draw (y)--(i);
      \draw [red!75!black] (r)--(i);
      \draw [green!75!black] (b)--(i);
      \draw [blue!75!black] (g)--(i);
      
      \draw (x)--(k);
      \draw [blue!75!black] (r)--(k);
      \draw [red!75!black] (b)--(k);
      \draw [green!75!black] (g)--(k);

      \draw (i)--(k);
    \end{tikzpicture}
    \caption{Analysis of the color of the edges incident on $X_{yrbg}$ and on
      $X_{xrbg}$.}\label{fig:rbgyx}
  \end{figure}
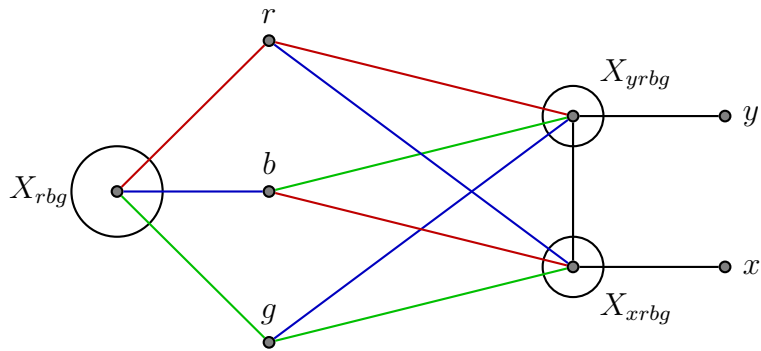

  Therefore, we have that each vertex in $X_{rbg} \cup N_G(X_{rbg})$ belongs to
  one of the monochromatic components $C_\red(r)$, $C_\blue(b)$ or
  $C_\green(g)$, while a vertex in $V(G)\setminus \left(X_{rbg} \cup
    N_G(X_{rbg}) \right)$ belongs to one of the monochromatic components
  $C_\red(r)$, $C_\green(b)$ or $C_\blue(g)$ where the latter two are the green
  component containing~$b$ and the blue component containing~$g$, respectively.
  This gives a covering of $G$ with five monochromatic trees. Next we will show
  that actually three of those trees already cover all the vertices.

  Suppose that at least~$4$ among the components
  $C_{\red}(r)$,~$C_{\blue}(b)$,~$C_{\green}(b)$,~$C_{\green}(g)$,
  and~$C_{\blue}(g)$ are needed to cover all vertices. Since there does not
  exist any monochromatic path between any two of~$r$, $b$, $g$, we know that
  for each~$i \in \{r,b,g\}$, any monochromatic component containing~$i$ does
  not intersect~$\{r,g,b\}\setminus\{i\}$. Hence, among those at least~$4$
  components, we have for each~$i\in\{r,b,g\}$ one component containing it and,
  without loss of generality, two containing~$b$. That is, three components of
  those at least $4$ components needed to cover all the vertices
  are~$C_{\red}(r)$,~$C_{\blue}(b)$ and~$C_{\green}(b)$. Now there are two cases
  regarding the fourth component: we need~$C_{\green}(g)$ as the fourth
  component or we need~$C_{\blue}(g)$ (those two cases might intersect).

  We begin with the first case, where we need the components $C_{\red}(r)$,
  $C_{\blue}(b)$, $C_{\green}(b)$ and $C_{\green}(g)$ to cover all the vertices
  of $G$. Let
  \[
    \tilde{b} \in C_{\blue}(b) \setminus \left( C_{\red}(r) \cup C_{\green}(b)
      \cup C_{\green}(g) \right)
  \]
  and let
  \[
    \tilde{g} \in C_{\green}(b)\setminus \left( C_{\red}(r) \cup C_{\blue}(b)
      \cup C_{green}(g) \right).
  \]
  Then let~$X_{\tilde{b}\tilde{g}rbg}$ be the maximum set of common neighbours
  of~$\tilde{b},\tilde{g},r,g,b$ such that for each~$i \in \{\tilde{b},
  \tilde{g}, r, b, g\}$, the edges from~$i$ to~$X_{\tilde{b}\tilde{g}rbg}$ are
  all coloured the same. Since~$\vert X_{\tilde{b}\tilde{g}rbg}\vert \geq
  np^5/240 \gg \frac{\log{n}}{p} $, we have
  \[
    |E_G(X_{\tilde{b}\tilde{g}rbg},X_{yrbg})| \geq 1 \qand
    |E_G(X_{\tilde{b}\tilde{g}rbg},X_{rbg})| \geq 1.
  \]
  We will analyse the possible colours of the edges between the specified
  vertices and~$X_{\tilde{b}\tilde{g}rbg}$. If for each of~${r,b,g}$, the colour
  it sends to~$X_{\tilde{b}\tilde{g}rbg}$ is different from the colour it sends
  to~$X_{rbg}$, then any edge between~$X_{\tilde{b}\tilde{g}rbg}$ and~$X_{rbg}$
  ensures a monochromatic path between two of~$r$,~$b$,~$g$ (in the colour of
  that edge). Similarly, it cannot happen that for each of~${r,b,g}$, the colour it sends
  to~$X_{\tilde{b}\tilde{g}rbg}$ is different from the colour it sends
  to~$X_{yrbg}$. Thus, since~$r$ sends red to both~$X_{rbg}$ and~$X_{yrbg}$ while the colours from~$b$ (and~$g$) to~$X_{rbg}$ and~$X_{yrbg}$ are switched, the colour of the edges between~$r$
  and~$X_{\tilde{b}\tilde{g}rbg}$ is red.
  
  Now note that, by the choice of~$\tilde{b}$ and~$\tilde{g}$, the edges between
  each of them and~$X_{\tilde{b}\tilde{g}rbg}$ can not be red. Further, the
  choice implies that an edge between~$\tilde{b}$
  and~$X_{\tilde{b}\tilde{g}rbg}$ can not be of the same colour (green or blue)
  as an edge between~$\tilde{g}$ and~$X_{\tilde{b}\tilde{g}rbg}$. If~$g$ would
  send blue (and hence~$b$ would send green) edges
  to~$X_{\tilde{b}\tilde{g}rbg}$, there would either be a blue path between~$b$
  and~$g$ (if the edges between~$\tilde{b}$ and~$X_{\tilde{b}\tilde{g}rbg}$ are
  blue) or~$\tilde{b}$ would lie in~$C_\green(b)$ (if the edges
  between~$\tilde{b}$ and~$X_{\tilde{b}\tilde{g}rbg}$ are green). Since both
  those situations would mean a contradiction, we may assume that each
  of~$r$,~$b$,~$g$ sends edges with that colour to~$X_{\tilde{b}\tilde{g}rbg}$
  as it does to~$X_{rbg}$. But then~$X_{\tilde{b}\tilde{g}rbg}$ is actually a
  subset of~$X_{rbg}$ and therefore~$\tilde{g}$, having an edge to~$X_{rbg}$,
  lies in one of~$C_{\red}(r)$, $C_\blue(b)$, or $C_\green(g)$, a contradiction.
  
  In the case where the forth component that we need is~$C_\blue(g)$, we repeat
  the construction of~$X_{\tilde{b}\tilde{g}rbg}$ similarly as before by letting
  \[
    \tilde{b} \in C_\blue(b) \setminus ( C_\red(r)\cup C_\green(b)\cup
    C_\blue(g))
  \]
  and
  \[
    \tilde{g}\in C_\green(b) \setminus ( C_\red(r)\cup C_\blue(b)\cup
    C_\blue(g)) .
  \]
  Also as before, we end up with~$X_{\tilde{b}\tilde{g}rbg}$ being part
  of~$X_{rbg}$. From the choice of~$\tilde{g}$, the edges it sends
  to~$X_{\tilde{b}\tilde{g}rbg}$ have to be green, since otherwise it would be
  in~$C_{\red}(r)$ or~$C_\blue(b)$. But that gives a green path between~$b$
  and~$g$, a contradiction.
		
  Summarising, we infer that three components
  among~$C_\red(r)$,~$C_\blue(b)$,~$C_\green(b)$,~$C_\green(g)$ and~$C_\blue(g)$
  cover the vertex set of~$G$.
\end{proof}

\subsection{Shortcut graphs with independence number at most two}\label{sec:alpha3}
\begin{proof}[Proof of~\cref{lemma:main} for~$\alpha(F)\leq 2$]
  We start by noticing that if~$\alpha (F)=1$, then the graph~$F$ together with
  the colouring~$\varphi'$ is a complete~$3$-coloured graph and therefore,
  by~\cref{lemma:EGP91}, there exists a partition of~$V(F)$ into~$2$
  monochromatic trees. Thus, we may assume that~$\alpha(F)=2$.

  Let $\mathcal{H}$ be the 3-uniform hypergraph with~$V(\mathcal{H})$ being the
  collection of all the monochromatic components of~$F$ under the
  colouring~$\varphi'$ and three monochromatic components form a hyperedge
  in~$\mathcal{H}$ if they share a vertex. Notice that~$\mathcal{H}$ is
  3-partite, since distinct monochromatic components of the same colour do not
  have a common vertex and therefore they can not belong to the same hyperedge.
  In other words, the colour of each component give us a 3-partition of the
  vertex set of~$\mathcal{H}$. We denote by~$V_{\red}$,$V_{\blue}$ and~$V_{\green}$ the set of vertices of~$V(\mathcal{H})$ that correspond to,
  respectively, red, blue and green components. Such construction was inspired
  by a construction due to Gy\'arf\'as~\cite{Gyarfas}.

  Note that every vertex~$v$ of~$F$ is contained in a monochromatic component
  for each one of the colours (a monochromatic component could consist
  only of $v$). Therefore, any vertex cover of~$\mathcal{H}$ corresponds to a
  covering of the vertices of~$F$ with monochromatic trees. Indeed, if~$A$ is a
  vertex cover of~$\mathcal{H}$, then consider the monochromatic components
  corresponding to each vertex in~$A$. If any vertex~$v$ of~$F$ is not covered
  by those components, then the vertices in $\mathcal{H}$ corresponding to the
  red, green and blue components in $F$ containing $v$ do not belong to $A$ and
  they form an hyperedge. But this contradicts the fact that $A$ is a vertex
  cover of $\mathcal{H}$. Therefore,
  \begin{align}\label{tctau}
    \tc(F,\phi') \leq \tau(\mathcal{H}).
  \end{align}

  Let~$L = \bigcup_{s\in V_{\red}}L_s$ be the union of the link graphs~$L_s$ of all
  vertices~$s\in V_{\red}$. Any vertex cover of this bipartite graph~$L$
  corresponds to a vertex cover of~$\mathcal{H}$ of the same size. Therefore,
  \begin{align}\label{tautau}
    \tau(\mathcal{H}) \leq \tau(L).
  \end{align}
  Furthermore, by K\"{o}nig's theorem we know that $\tau(L) = \nu(L)$. Thus, if
  $\nu(L) \leq 3$, then by~\eqref{tctau} and~\eqref{tautau}, we have
  \begin{align*}
    \tc(F,\phi') \leq \tau(\mathcal{H}) \leq \tau(L) = \nu(L) \leq 3.
  \end{align*}
  Therefore, we may assume that \(\nu(L) \geq 4\), and fix a matching~$M_L$ of
  size at least~$4$ in~$L$. Let us say that $M_L$ consists of the edges
  $G_1B_1$, $G_2B_2$, $G_3B_3$, and $G_4B_4$, where $\{G_1,G_2,G_3,G_4\}
  \subseteq V_{\green}$ and $\{B_1,B_2,B_3,B_4\} \subseteq V_{\blue}$.

  Now we give an upper bound for~$\nu(\mathcal{H})$. Note that any
  matching~$M_{\mathcal{H}}$ in~$\mathcal{H}$ gives us an independent set~$I$ in~$F$.
  Indeed, for each hyperedge~$e \in M_{\mathcal{H}}$, let~$v_e\in V(F)$ be any vertex in the
  intersection of those monochromatic components associated to the vertices
  in~$e$ and let~$I = \{v_e : e \in M_{\mathcal{H}}\}$. We claim that~$I$ is an independent
  set in~$F$. Indeed, if~$v_e$ and $v_f$ were adjacent vertices in~$I$, then~$e$
  and~$f$ intersect, as the edge connecting~$v_e$ to~$v_f$
  in~$F$ will connect the monochromatic components containing~$v_e$ and~$v_f$ of
  that colour that is given to the edge~$v_e v_f$. Therefore,
  since~$\alpha(F) = 2$, we have
  \begin{align}\label{nual}
    \nu(\mathcal{H}) \leq \alpha(F) = 2.
  \end{align}
  Now, if there are three different edges in~$M_L$ that are edges in the link
  graphs of three different vertices of~$V_{\red}$, then there would be a
  matching of size~$3$ in~$\mathcal{H}$, contradicting~\eqref{nual}. Therefore,
  we may assume that~$M_L$ is contained in the union of at most two link graphs,
  say~$L_{R_1}$ and~$L_{R_2}$, of vertices $R_1,R_2 \in V_\red$. Now we are left
  with three cases: (Case~\ref{case_1}) two edges of~$M_L$ belong to~$L_{R_1}$
  and two belong to~$L_{R_2}$; (Case~\ref{case_2}) three edges of~$M_L$ belong
  to~$L_{R_1}$ and one to~$L_{R_2}$; (Case~\ref{case_3}) the four edges of $M_L$
  belong to~$L_{R_{1}}$. Without loss of generality, we can describe each of
  those three cases as follows (see
  Figures~\ref{fig:case1},~\ref{fig:case2} and~\ref{fig:case3}):

  \begin{case}\label{case_1}
    The edges $G_1B_1$ and $G_2B_2$ belong to~$L_{R_1}$ and the edges $G_3B_3$
    and $G_4B_4$ belong to~$L_{R_2}$. That means that all the following four
    sets are non-empty:
    \begin{align*}
      J_1&:=R_1\cap G_1\cap B_1,\\
      J_2&:=R_1\cap G_2\cap B_2,\\
      J_3&:=R_2\cap G_3\cap B_3,\\
      J_4&:=R_2\cap G_4\cap B_4.
    \end{align*}
  \end{case}

  \begin{case}\label{case_2}
    The edges $G_1B_1$, $G_2B_2$ and $G_3B_3$ belong to~$L_{R_1}$ and the edge
    $G_4B_4$ belongs to~$L_{R_2}$. That means that all the following four sets
    are non-empty:
    \begin{align*}
      J_1&:=R_1\cap G_1\cap B_1,\\
      J_2&:=R_1\cap G_2\cap B_2,\\
      J_3&:=R_1\cap G_3\cap B_3,\\
      J_4&:=R_2\cap G_4\cap B_4.
    \end{align*}
  \end{case}

  \begin{case}\label{case_3}
    The edges $G_1B_1$, $G_2B_2$, $G_3B_3$ and $G_4B_4$ belong to~$L_{R_1}$.
    That means that all the following four sets are non-empty:
    \begin{align*}
      J_1&:=R_1\cap G_1\cap B_1,\\
      J_2&:=R_1\cap G_2\cap B_2,\\
      J_3&:=R_1\cap G_3\cap B_3,\\
      J_4&:=R_1\cap G_4\cap B_4.
    \end{align*}

    In this case, let~$R_2$ be any other red component different from~$R_1$ and
    let~$B$ and~$G$ be, respectively, a blue and a green component with~$R_2\cap
    B \cap G\neq \emptyset$. Suppose that $G \notin \{G_1,G_2,G_3,G_4\}$. Then
    the three of the edges $G_1,B_1$, $G_2,B_2$, $G_3,B_3$ and $G_4,B_4$ are not
    incident to $GB$ (because $B$ must be different of at least three of the
    sets $B_1$, $B_2$, $B_3$ and $B_4$) and those three edges together with $GB$
    may be analysed just as in Case~\ref{case_2}. Therefore, we may suppose that
    $G \in \{G_1,G_2,G_3,G_4\}$. Let us say, without loss of generality, that $G
    = G_4$. If $B \notin\{B_1,B_2,B_3\}$, then the edges $G_1B_1$, $G_2B_2$ and
    $G_3B_3$ belong to $L_{R_1}$, the edge $GB$ belongs to $L_{R_2}$ and this
    case may be analysed, again, just as in Case~\ref{case_2}. Therefore, we may
    assume that $B \in \{B_1,B_2,B_3\}$. Let us say, without loss of generality
    that $B = B_3$. Then let $J_5$ be the following non-empty set:
    \begin{align}\label{eq:case3-extra}
      J_5:=R_2\cap G_4\cap B_3.
    \end{align}
  \end{case}

  Let us further remark that, since~$\nu (\mathcal{H})\leq 2$, in each of the
  three cases above, we have
  \begin{align*}
    V(F) = R_1 \cup R_2 \cup G_1 \cup G_2 \cup G_3 \cup G_4 \cup B_1 \cup B_2
    \cup B_3 \cup B_4.
  \end{align*}
  Otherwise, for any uncovered vertex~$v \in V(F)$, the hyperedge given by the
  red, blue and green components containing $v$ together with the hyperedges
  $R_1B_1G_1$ and $R_2B_3G_3$ (in Cases~\ref{case_1} and~\ref{case_2}) or~$R_2B_3G_4$ (in Case~\ref{case_3}) give a matching of size~$3$ in~$\mathcal{H}$.

  \begin{figure}
    \centering
    \begin{tikzpicture}[scale=1, thick, auto, vertex/.style={circle, draw,
        fill=black!50, inner sep=0pt, minimum width=4pt}]
      
      \node [label={[label distance=.8cm]left:$R_1$}](r1) at (0,0) {};
      \node [label={[label distance=.8cm]right:$R_2$}](r2) at (4,0) {};
      
      \node [label={[label distance=.1cm]-45:$B_1$}](b1) at (0.25,-1) {};
      \node [label={[label distance=.1cm]-135:$G_1$}](g1) at (-0.25,-1) {};
      \node [vertex, label={[label distance=.1cm]above:$j_1$}](j1) at (0,-0.85) {};
      
      \node [label={[label distance=.1cm]45:$B_2$}](b2) at (0.25,1) {};
      \node [label={[label distance=.1cm]135:$G_2$}](g2) at (-0.25,1) {};
      \node [vertex, label={[label distance=.1cm]below:$j_2$}](j2) at (0,0.85) {};
      
      \node [label={[label distance=.1cm]-135:$B_3$}](b3) at (3.75,-1) {};
      \node [label={[label distance=.1cm]-45:$G_3$}](g3) at (4.25,-1) {};
      \node [vertex, label={[label distance=.1cm]above:$j_3$}](j3) at (4,-0.85) {};
      
      \node [label={[label distance=.1cm]135:$B_4$}](b4) at (3.75,1) {};
      \node [label={[label distance=.1cm]45:$G_4$}](g4) at (4.25,1) {};
      \node [vertex, label={[label distance=.1cm]below:$j_4$}](j4) at (4,0.85) {};

      \draw[red!75!black, line width=1pt] (r1) circle (1 cm);
      \draw[red!75!black, line width=1pt] (r2) circle (1 cm);
      
      \draw[blue!75!black, line width=1pt] (b1) circle (0.4 cm);
      \draw[green!75!black, line width=1pt] (g1) circle (0.4 cm);

      \draw[blue!75!black, line width=1pt] (b2) circle (0.4 cm);
      \draw[green!75!black, line width=1pt] (g2) circle (0.4 cm);
      
      \draw[blue!75!black, line width=1pt] (b3) circle (0.4 cm);
      \draw[green!75!black, line width=1pt] (g3) circle (0.4 cm);
      
      \draw[blue!75!black, line width=1pt] (b4) circle (0.4 cm);
      \draw[green!75!black, line width=1pt] (g4) circle (0.4 cm);
      
    \end{tikzpicture}
    \caption{Case 1}\label{fig:case1}
  \end{figure}

  Let us start with Case~\ref{case_1}.

\medskip

\noindent\textit{Proof in Case~\ref{case_1}}:
  We will prove that~$R_1$ and~$R_2$ together with possibly one further
  monochromatic component cover $V(F)$. For each $i \in \{1,2,3,4\}$, let $\tilde{B}_i =
  B_i\setminus(R_1\cup R_2)$ and $\tilde{G}_i = G_i\setminus (R_1 \cup R_2)$.

  Pick vertices~$j_i\in J_i$, with~$i\in\{1,2,3,4\}$, arbitrarily. Consider a
  vertex~$o\in \tilde{B}_1$ (if such a vertex exists). Since~$\alpha (F) = 2$,
  there is an edge connecting two of~$o$,~$j_2$,~$j_3$. Because $j_2$ and $j_3$
  belong to different components of each colour, such an edge must be incident to
  $o$. So let us say that such edge is $oj_i$, for some $i\in\{2,3\}$. Since~$o
  \notin R_1 \cup R_2$, the edge $oj_i$ cannot be red. And since $o \in B_1$,
  $oj_i$ cannot be blue either, otherwise we would connect the blue
  components $B_1$ and $B_i$. Now assume that $o$ and $j_2$ are not adjacent.
  Then $oj_3$ is a green edge in $F$. By analogously analysing the edge between
  $o$, $j_2$ and $j_4$ together with the supposition that $oj_2$ is not an edge
  in $F$, we get that $oj_4$ must be a green edge in $F$. But then we have a
  green path $j_3oj_4$ connecting $j_3$ to $j_4$, a contradiction. Therefore
  $oj_2$ is an edge in $F$ and it is green. That implies that $o \in G_2$.
  Therefore $\tilde{B}_1 \subseteq G_2$. Analogously, we can conclude the
  following:
  \begin{equation}
    \begin{aligned}\label{conn_comp}
      & \tilde{B}_1 \subseteq G_2, & \tilde{G}_1 \subseteq B_2, \\
      & \tilde{B}_2 \subseteq G_1, & \tilde{G}_2 \subseteq B_1, \\
      & \tilde{B}_3 \subseteq G_4, & \tilde{G}_3 \subseteq B_4, \\
      & \tilde{B}_4 \subseteq G_3, & \tilde{G}_4 \subseteq B_3. \\
    \end{aligned}
  \end{equation}

  \begin{claim}
    We have $\tilde{B}_1 \cup \tilde{G}_1 \cup \tilde{B}_2 \cup \tilde{G}_2 =
    \emptyset$ or $\tilde{B}_3 \cup \tilde{G}_3 \cup \tilde{B}_4 \cup
    \tilde{G}_4 = \emptyset$.
  \end{claim}
  \begin{claimproof}
    Suppose for a contradiction that there exist~$o_1\in \tilde{B}_1 \cup
    \tilde{G}_1 \cup \tilde{B}_2 \cup \tilde{G}_2$ and~$o_2\in \tilde{B}_3 \cup
    \tilde{G}_3 \cup \tilde{B}_4 \cup \tilde{G}_4$. Recall that from our choice
    of $p$, there is some~$z\in N(j_1,j_2,j_3,j_4,o_1,o_2)$. Two of the edges
    $zj_i$,for $i \in \{1,2,3,4\}$, have the same colour. Since each $j_i$ belongs to
    different green and blue components, those two edges are red. Since
    $\{j_1,j_2\} \in R_1$ and $\{j_3,j_4\}\in R_2$, those two red edges are
    either $zj_1$ and $zj_2$ or $zj_3$ and $zj_4$. Let us say that $zj_1$ and
    $zj_2$ are red (the other case is similar). Then one of the edges $zj_3$ and
    $zj_4$ has to be green and the other blue. Now, since~$o_1 \notin R_1$, the
    edge $zo_1$ is either green or blue. Then one of the paths~$o_1zj_3$ or
    $o_1zj_4$ is green or blue. This implies that $o_1 \in B_3\cup G_3\cup
    B_4\cup G_4$. On the other hand,~\eqref{conn_comp} implies that $o_1 \in
    \left(B_1\cup B_2\right) \cap \left(G_1 \cup G_2\right)$. But then we
    reached a contradiction, since that would mean that $o_1$ belongs to two
    different components of the same colour.
  \end{claimproof}

  We may assume without loss of generality that~$\tilde{B}_3 \cup \tilde{G}_3
  \cup \tilde{B}_4 \cup \tilde{G}_4$ is empty. Then, recalling
  that~$\nu(\mathcal{H})\leq 2$ and in view of~\eqref{conn_comp}, the union of
  the components~$R_1$, $B_1$, $G_1$ and $R_2$ covers every vertex of $F$. If we
  show that ${B}_1\subseteq G_1 \cup R_1 \cup R_2$ or that ${G}_1\subseteq B_1
  \cup R_1 \cup R_2$, then we get three monochromatic components covering the
  vertices of $F$. Our next claim states precisely that.

  \begin{claim}
    We have $\tilde{B}_1\setminus G_1 = \emptyset$ or $\tilde{G}_1 \setminus
    B_1= \emptyset$.
  \end{claim}
  \begin{claimproof}
    Suppose that there exist two distinct vertices~$b\in \tilde{B}_1\setminus
    G_1$ and~$g\in \tilde{G}_1\setminus B_1$. Let~$z\in N(j_1,j_2,j_3,j_4,b,g)$.
    As before, either~$zj_1$ and~$zj_2$ or~$zj_3$ and~$zj_4$ are red edges. First assume that $zj_1$ and $zj_2$ are red. Then one of the edges $zj_3$ and
    $zj_4$ has to be green and the other blue. Now, since $b \notin R_1$, the
    edge $zb$ is either green or blue. Then one of the paths $bzj_3$ or $bzj_4$
    is green or blue. This implies that $b \in B_3\cup G_3\cup B_4\cup G_4$. On
    the other hand,~\eqref{conn_comp} implies that $b \in B_1 \cap G_2$. Then we
    reached a contradiction, since that would mean that $b$ belongs to two
    different components of the same colour.

    Therefore, the edges $zj_3$ and $zj_4$ are red and one of the edges $zj_1$
    and $zj_2$ is green and the other is blue. First let us say that $zj_1$ is
    green and $zj_2$ is blue. Since $b\notin (R_1\cup R_2)$, the edge $zb$
    cannot be red. Also the edge $zb$ cannot be blue otherwise the path $bzj_2$
    would connect the components $B_1$ and $B_2$. Finally, $zb$ cannot be green,
    otherwise the path~$bzj_1$ would gives us that $b\in G_1$. Therefore $zj_1$
    is blue and $zj_2$ is green. But this case analogously leads to a
    contradiction (with~$g$ and~$G_i$ instead of~$b$ and~$B_i$ and green and
    blue switched).
  \end{claimproof}

  \begin{figure}
    \centering
    \begin{tikzpicture}[scale=1, thick, auto, vertex/.style={circle, draw,
        fill=black!50, inner sep=0pt, minimum width=4pt}]
      
      \node [label={[label distance=.8cm]left:$R_1$}](r1) at (0,0) {};
      \node [label={[label distance=.8cm]right:$R_2$}](r2) at (4,0) {};
      
      \node [label={[label distance=.1cm]-45:$B_1$}](b1) at (0.25,-1) {};
      \node [label={[label distance=.1cm]-135:$G_1$}](g1) at (-0.25,-1) {};
      \node [vertex, label={[label distance=.1cm]above:$j_1$}](j1) at (0,-0.85) {};
      
      \node [label={[label distance=.1cm]45:$B_2$}](b2) at (0.25,1) {};
      \node [label={[label distance=.1cm]135:$G_2$}](g2) at (-0.25,1) {};
      \node [vertex, label={[label distance=.1cm]below:$j_2$}](j2) at (0,0.85) {};
      
      \node [label={[label distance=.1cm]15:$B_3$}](b3) at (1,0.25) {};
      \node [label={[label distance=.1cm]-15:$G_3$}](g3) at (1,-0.25) {};
      \node [vertex, label={[label distance=.1cm]left:$j_3$}](j3) at (0.9,0) {};
      
      \node [label={[label distance=.1cm]135:$B_4$}](b4) at (3.75,1) {};
      \node [label={[label distance=.1cm]45:$G_4$}](g4) at (4.25,1) {};
      \node [vertex, label={[label distance=.1cm]below:$j_4$}](j4) at (4,0.85) {};

      \draw[red!75!black, line width=1pt] (r1) circle (1 cm);
      \draw[red!75!black, line width=1pt] (r2) circle (1 cm);
      
      \draw[blue!75!black, line width=1pt] (b1) circle (0.4 cm);
      \draw[green!75!black, line width=1pt] (g1) circle (0.4 cm);
      
      \draw[blue!75!black, line width=1pt] (b2) circle (0.4 cm);
      \draw[green!75!black, line width=1pt] (g2) circle (0.4 cm);
      
      \draw[blue!75!black, line width=1pt] (b3) circle (0.4 cm);
      \draw[green!75!black, line width=1pt] (g3) circle (0.4 cm);
      
      \draw[blue!75!black, line width=1pt] (b4) circle (0.4 cm);
      \draw[green!75!black, line width=1pt] (g4) circle (0.4 cm);
      
    \end{tikzpicture}
    \caption{Case 2}\label{fig:case2}
  \end{figure}

We proceed to the proof of Case~\ref{case_2}.

\medskip
		
\noindent\textit{Proof in Case~\ref{case_2}}:
  As in Case~\ref{case_1}, pick vertices~$j_i\in J_i$, with~$i\in\{1,2,3,4\}$
  arbitrarily. We claim that~$V(F)\subseteq R_1 \cup R_2 \cup B_4 \cup G_4$.
  Indeed, let~$o\in V(F)\setminus (R_1\cup R_2)$. Notice that
  since~$\alpha(F)=2$, there is an edge in each of the following sets of three
  vertices:~$\{o, j_4, j_1\}$,~$\{o, j_4, j_2\}$, and~$\{o, j_4, j_3\}$. We
  claim that~$oj_4$ is an edge of $F$. Indeed, if this was not the case, then
  since there cannot be an edge between~$j_4$ and~$j_i$ for~$i=1,2,3$, we would
  have the edges~$oj_1$,~$oj_2$ and~$oj_3$ and all of them would be coloured
  green or blue. Thus, two of them would be coloured the same, connecting two
  distinct components of one colour in this colour, a contradiction. So~$oj_4\in
  E(F)$ and since~$oj_4$ cannot be red, we conclude that~$o\in (B_4\cup G_4)$.
  Therefore,~$R_1$,~$R_2$,~$B_4$ and~$G_4$ cover all vertices of~$F$.

  If $B_4\setminus (R_1\cup R_2\cup G_4)=\emptyset$ or $G_4\setminus (R_1\cup
  R_2\cup B_4)=\emptyset$, then we get three monochromatic components covering
  $V(F)$. So let us assume that there exist~$b\in B_4\setminus (R_1\cup R_2\cup
  G_4)$ and~$g\in G_4\setminus (R_1\cup R_2\cup B_4)$. If $b$ and~$g$ are not
  adjacent, then since each of the sets~$\{b, g, j_i\}$, for~$i=1,2,3$, has to
  induce at least one edge, there are two edges between $b$
  and~$\{j_1,j_2,j_3\}$ or two edges between $g$ and~$\{j_1,j_2,j_3\}$. However,
  from the choice of $b$, we know that all the edges between $b$
  and~$\{j_1,j_2,j_3\}$ are green, and therefore two of such edges would give us
  a green connection between two different green components, a contradiction.
  Similarly, from the choice of $g$, we know that all the edges between $b$ and
  $\{j_1,j_2,j_3\}$ are blue, and two of such edges would give us a blue
  connection between two different blue components, again a contradiction.

  Hence, we conclude that $bg\in F$ for any $b\in B_4\setminus (R_1\cup R_2\cup
  G_4)$ and any~$g\in G_4\setminus (R_1\cup R_2\cup B_4)$ and any such edge $bg$
  is red. Therefore, there is a red component~$R_3$ covering~$(B_4\triangle
  G_4)\setminus (R_1\cup R_2)$, where~$B_4\triangle G_4$ denotes the symmetric
  difference. If~$(B_4\cap G_4) \setminus (R_1\cup R_2) =\emptyset$,
  then~$R_1$,~$R_2$ and~$R_3$ cover~$V(F)$ and we are done. Therefore, suppose
  there is a vertex~$x\in (B_4\cap G_4) \setminus (R_1\cup R_2)$.
  If~$R_2\setminus (B_4\cup G_4)=\emptyset$, then~$R_1$,~$B_4$,~$G_4$
  cover~$V(F)$ and we are done. Therefore, suppose there is a vertex $y\in
  R_2\setminus (B_4\cup G_4)$. Note that~$xy\notin E(F)$, since $x$ and $y$
  belong to different components in each of the colours. Also, $xj_i\notin
  E(F)$, for~$i\in\{1,2,3\}$, since otherwise two different components of the
  same colour would be connected in that colour by the edge $xj_i$. Now~$\alpha
  (F)=2$ implies that~$yj_i\in E(F)$, for~$i\in\{1,2,3\}$ (otherwise,
  $\{x,y,j_i\}$ would be an independent set). But these edges must all be green
  or blue, hence two of them are of the same colour, connecting two different
  components of one colour in that colour, a contradiction.

  \begin{figure}
    \centering
    \begin{tikzpicture}[scale=1, thick, auto, vertex/.style={circle, draw,
        fill=black!50, inner sep=0pt, minimum width=4pt}]
      
      \node [label={[label distance=.8cm]left:$R_1$}](r1) at (0,0) {};
      \node [label={[label distance=.8cm]right:$R_2$}](r2) at (4,0) {};
      
      \node [label={[label distance=.2cm]below:$B_1$}](b1) at (-0.25,-0.9) {};
      \node [label={[label distance=.1cm]-135:$G_1$}](g1) at (-0.75,-0.7) {};
      \node [vertex, label={[label distance=.1cm]above:$j_1$}](j1) at (-0.5,-0.75) {};
      
      \node [label={[label distance=.2cm]above:$B_2$}](b2) at (-0.25,0.9) {};
      \node [label={[label distance=.1cm]135:$G_2$}](g2) at (-0.75,0.7) {};
      \node [vertex, label={[label distance=.1cm]below:$j_2$}](j2) at (-0.5,0.75) {};
      
      \node [label={[label distance=.1cm]below:$B_3$}](b3) at (2,-0.25) {};
      \node [label={[label distance=.1cm]below:$G_3$}](g3) at (0.75,-0.7) {};
      \node [vertex, label={[label distance=.1cm]left:$j_3$}](j3) at (0.7,-0.5) {};
      
      \node [label={[label distance=.1cm]above:$B_4$}](b4) at (0.75,0.7) {};
      \node [label={[label distance=.1cm]above:$G_4$}](g4) at (2,0.25) {};
      \node [vertex, label={[label distance=.1cm]left:$j_4$}](j4) at (0.7,0.5) {};
      
      \node [vertex, label={[label distance=.2cm]right:$j_5$}](j5) at (3.2,0) {};
      
      \draw[red!75!black, line width=1pt] (r1) circle (1 cm);
      \draw[red!75!black, line width=1pt] (r2) circle (1 cm);
      
      \draw[blue!75!black, line width=1pt] (b1) circle (0.4 cm);
      \draw[green!75!black, line width=1pt] (g1) circle (0.4 cm);
      
      \draw[blue!75!black, line width=1pt] (b2) circle (0.4 cm);
      \draw[green!75!black, line width=1pt] (g2) circle (0.4 cm);
      
      \draw[rotate around={-80:(b3)},blue!75!black, line width=1pt] (b3) ellipse (8pt and 1.6cm);
      \draw[green!75!black, line width=1pt] (g3) circle (0.4 cm);
      
      \draw[rotate around={80:(g4)},green!75!black, line width=1pt] (g4) ellipse (8pt and 1.6cm);
      \draw[blue!75!black, line width=1pt] (b4) circle (0.4 cm);
      
    \end{tikzpicture}
    \caption{Case 3}\label{fig:case3}
  \end{figure}

We arrived at the last case, Case~\ref{case_3}.

\medskip

\noindent\textit{Proof in Case~\ref{case_3}}:
  Similarly to the previous cases, let us pick vertices~$j_i\in J_i$,
  with~$i\in\{1,2,3,4,5\}$ arbitrarily. We will show first that we can cover all vertices
  of~$F$ with~$4$ monochromatic components. Let~$o_1,o_2 \in V(F)\setminus
  (R_1\cup B_3\cup G_4)$ and let~$z\in N(j_1,j_2,j_3,o_1,o_2,j_5)$. At least one
  of the edges~$zj_1$,~$zj_2$ and~$zj_3$ is red, as otherwise we would connect
  two distinct components of one colour in that colour. Therefore~$z \in R_1$. Since~$o_1,o_2,j_5 \notin R_1$, the edges~$zo_1$,~$zo_2$ and~$zj_5$ cannot be red.
  Furthermore,~$o_1z$ and~$o_2z$ are coloured with a colour different from the
  colour of the edge~$j_5z$, as otherwise they would belong to~$B_3$ or~$G_4$.
  Thus,~$o_1$ and~$o_2$ are connected by a monochromatic path in green or blue.
  Hence, we showed that any two vertices of~$V(F)\setminus (R_1\cup B_3\cup
  G_4)$ are connected by a monochromatic path in green or blue. We infer
  that there is a green or blue component covering~$V(F)\setminus (R_1\cup
  B_3\cup G_4)$. Therefore,~$R_1$,~$B_3$,~$G_4$ and one further blue or green
  component~$C$ cover all vertices of~$G$. Let us assume that~$C$ is a green
  component; the case where~$C$ is a blue component is analogous.

  We claim that~$R_1 \cup B_3 \cup C$, or~$R_1 \cup G_4 \cup C$, or~$R_1 \cup
  B_3 \cup G_4$ covers~$V(F)$. Indeed, suppose for the sake of contradiction that
  there exist vertices~$g\in G_4\setminus (R_1\cup B_3\cup C)$,~$b\in
  B_3\setminus (R_1\cup G_4\cup C)$ and~$c\in C\setminus (R_1\cup B_3\cup G_4)$.
  Let~$z\in N(j_1,j_2,j_3,g,b,c)$ and note that one of~$zj_1$,~$zj_2$ and~$zj_3$
  is red. Consequently~$gz$,~$cz$ and~$bz$ are not red. Notice, however, that~$gz$ and~$bz$ can not be both green and neither both blue. Now let us say~$cz$
  is green. Since~$c \notin G_4$ and~$g \in G_4$, we would have~$gz$ blue in
  this case. But then~$bz$ must be green and since~$c\in C$ and~$C$ is a green
  component, we have~$b \in C$, which is a contradiction. Therefore~$cz$
  must be blue. Then, since~$c \notin B_3$ and~$b \in B_3$, the edge~$bz$ should
  be green. Thus the edge~$gz$ is blue. Since this argument holds for any~$g\in
  G_4\setminus (R_1\cup B_3\cup C)$ and~$c\in C\setminus (R_1\cup B_3\cup G_4)$,
  we conclude that~$V(F)\setminus (R_1\cup B_3)$ can be covered by one blue
  tree. Hence,~$G$ can be covered by the three monochromatic trees.
This finishes the last case and thereby the proof of~\cref{lemma:main}.
\end{proof}

\bibliography{bib}
\end{document}

%% file: preamble.tex
\pdfoutput=1
\makeatletter
\let\origsection=\section \def\section{\@ifstar{\origsection*}{\mysection}}
\def\mysection{\@startsection{section}{1}\z@{.7\linespacing\@plus\linespacing}{.5\linespacing}{\normalfont\scshape\centering\S}}
\makeatother

\usepackage{amsmath,amssymb,amsthm}
\usepackage{mathrsfs}
\usepackage{mathabx}\changenotsign
\usepackage{dsfont}

\usepackage{xcolor}
\usepackage[backref]{hyperref}
\hypersetup{
    colorlinks,
    linkcolor={red!60!black},
    citecolor={green!60!black},
    urlcolor={blue!60!black}
}

\usepackage[capitalise,noabbrev]{cleveref}

\usepackage[open,openlevel=2,atend]{bookmark}

\usepackage[abbrev,msc-links,backrefs]{amsrefs}
\usepackage{doi}

\renewcommand{\PrintDOI}[1]{\doi{#1}}

\usepackage[T1]{fontenc}
\usepackage{lmodern}

\usepackage{geometry}
\numberwithin{equation}{section}
\numberwithin{figure}{section}

\usepackage{enumitem}
\def\rmlabel{\upshape({\itshape \roman*\,})}

\let\log=\ln
\let\polishlcross=\l
\def\l{\ifmmode\ell\else\polishlcross\fi}

\def\qand{\quad\text{and}\quad}

\let\emptyset=\varnothing
\let\setminus=\smallsetminus

\makeatletter
\def\moverlay{\mathpalette\mov@rlay}
\def\mov@rlay#1#2{\leavevmode\vtop{   \baselineskip\z@skip \lineskiplimit-\maxdimen
   \ialign{\hfil$\m@th#1##$\hfil\cr#2\crcr}}}
\newcommand{\charfusion}[3][\mathord]{
    #1{\ifx#1\mathop\vphantom{#2}\fi
        \mathpalette\mov@rlay{#2\cr#3}
      }
    \ifx#1\mathop\expandafter\displaylimits\fi}
\makeatother

\DeclareFontFamily{U}  {MnSymbolC}{}
\DeclareSymbolFont{MnSyC}         {U}  {MnSymbolC}{m}{n}
\DeclareFontShape{U}{MnSymbolC}{m}{n}{
    <-6>  MnSymbolC5
   <6-7>  MnSymbolC6
   <7-8>  MnSymbolC7
   <8-9>  MnSymbolC8
   <9-10> MnSymbolC9
  <10-12> MnSymbolC10
  <12->   MnSymbolC12}{}
\DeclareMathSymbol{\powerset}{\mathord}{MnSyC}{180}

\usepackage{tikz}
\usetikzlibrary{calc,decorations.pathmorphing}
\pgfdeclarelayer{background}
\pgfdeclarelayer{foreground}
\pgfdeclarelayer{front}
\pgfsetlayers{background,main,foreground,front}

\usepackage{subcaption}
\newcommand{\qedge}[7]{

	\ifx\relax#4\relax
		\def\qoffs{0pt}
	\else
		\def\qoffs{#4}
	\fi

	\def\qhedge{
		($#1+#3!\qoffs!-90:#2-#3$) --
		($#2+#1!\qoffs!-90:#3-#1$) --
		($#3+#2!\qoffs!-90:#1-#2$) -- cycle}

	\coordinate (12) at ($#1!\qoffs!90:#2$);
	\coordinate (13) at ($#1!\qoffs!-90:#3$);
	\coordinate (23) at ($#2!\qoffs!90:#3$);
	\coordinate (21) at ($#2!\qoffs!-90:#1$);
	\coordinate (31) at ($#3!\qoffs!90:#1$);
	\coordinate (32) at ($#3!\qoffs!-90:#2$);
	
	\def\nqhedge{
		(13) let \p1=($(13)-#1$), \p2=($(12)-#1$) in
			arc[start angle={atan2(\y1,\x1)}, delta angle={atan2(\y2,\x2)-atan2(\y1,\x1)-360*(atan2(\y2,\x2)-atan2(\y1,\x1)>0)}, x radius=\qoffs, y radius=\qoffs] --
		(21) let \p1=($(21)-#2$), \p2=($(23)-#2$) in
			arc[start angle={atan2(\y1,\x1)}, delta angle={atan2(\y2,\x2)-atan2(\y1,\x1)-360*(atan2(\y2,\x2)-atan2(\y1,\x1)>0)}, x radius=\qoffs, y radius=\qoffs] --
		(32) let \p1=($(32)-#3$), \p2=($(31)-#3$) in
			arc[start angle={atan2(\y1,\x1)}, delta angle={atan2(\y2,\x2)-atan2(\y1,\x1)-360*(atan2(\y2,\x2)-atan2(\y1,\x1)>0)}, x radius=\qoffs, y radius=\qoffs] --
		cycle}

		\ifx\relax#5\relax
		\def\qlwidth{1pt}
	\else
		\def\qlwidth{#5}
	\fi
	
		\ifx\relax#7\relax
		\fill \nqhedge;
	\else
		\fill[#7]\nqhedge;
	\fi

		\ifx\relax#6\relax
		\draw[line width=\qlwidth,rounded corners=\qoffs]\nqhedge;
	\else
		\draw[line width=\qlwidth,#6]\nqhedge;
	\fi
}

\def\red{\text{\rm red}}
\def\blue{\text{\rm blue}}
\def\green{\text{\rm green}}
\def\tc{\text{\rm tc}}
\def\tp{\text{\rm tp}}

\let\epsilon=\varepsilon
\let\eps=\epsilon
\let\phi=\varphi
\let\rho=\varrho
\let\theta=\vartheta

\newtheoremstyle{note}  {4pt}  {4pt}  {\sl}  {}  {\bfseries}  {.}  {.5em}          {}
\newtheoremstyle{introthms}  {3pt}  {3pt}  {\itshape}  {}  {\bfseries}  {.}  {.5em}          {\thmnote{#3}}
\newtheoremstyle{remark}  {2pt}  {2pt}  {\rm}  {}  {\bfseries}  {.}  {.3em}          {}

\theoremstyle{plain}
\newtheorem{theorem}{Theorem}[section]
\newtheorem{lemma}[theorem]{Lemma}

\newtheorem{claim}[theorem]{Claim}

\theoremstyle{note}

\newtheorem{definition}[theorem]{Definition}

\theoremstyle{remark}

\usepackage{accents}

\usepackage{lineno}
\newcommand*\patchAmsMathEnvironmentForLineno[1]{
\expandafter\let\csname old#1\expandafter\endcsname\csname #1\endcsname
\expandafter\let\csname oldend#1\expandafter\endcsname\csname end#1\endcsname
\renewenvironment{#1}
{\linenomath\csname old#1\endcsname}
{\csname oldend#1\endcsname\endlinenomath}}
\newcommand*\patchBothAmsMathEnvironmentsForLineno[1]{
\patchAmsMathEnvironmentForLineno{#1}
\patchAmsMathEnvironmentForLineno{#1*}}
\AtBeginDocument{
\patchBothAmsMathEnvironmentsForLineno{equation}
\patchBothAmsMathEnvironmentsForLineno{align}
\patchBothAmsMathEnvironmentsForLineno{flalign}
\patchBothAmsMathEnvironmentsForLineno{alignat}
\patchBothAmsMathEnvironmentsForLineno{gather}
\patchBothAmsMathEnvironmentsForLineno{multline}
}

\newcommand{\gnp}{G(n,p)}

\newcommand{\oldqed}{}
\def\endofFact{\hfill\scalebox{.6}{$\Box$}}
\newenvironment{claimproof}[1][Proof]{
  \renewcommand{\oldqed}{\qedsymbol}
  \renewcommand{\qedsymbol}{\endofFact}
  \begin{proof}[#1]
}{
  \end{proof}
  \renewcommand{\qedsymbol}{\oldqed}
}

\def\qand{\quad\text{and}\quad}